\documentclass[12pt,a4paper]{article}
\usepackage{epsfig} 
\usepackage{float}
\usepackage{graphicx}
\usepackage{amsmath}
\usepackage{amssymb}
\usepackage{setspace}
\usepackage{lineno}

\usepackage{subcaption}

\usepackage{amsthm}
\usepackage{hyperref}

\theoremstyle{plain}
\newtheorem{theorem}{Theorem}[section]

\theoremstyle{definition}

\newtheorem{example}[theorem]{Example}

\theoremstyle{remark}
\newtheorem{remark}[theorem]{Remark}

\usepackage{xcolor}
\usepackage{soul}
\definecolor{akgcol}{rgb}{0,0,0.8}
\definecolor{afcol}{rgb}{1,0,0}

\begin{document}
\title{Diffusion on middle-$\xi$ Cantor sets}
\author{Alireza Khalili Golmankhaneh $^{a*}$, Arran Fernandez $^{b\dag}$\\ Ali Khalili Golmankhaneh $^{c\ddag}$ Dumitru Baleanu $^{d,f}$\S }
\date{}
\maketitle \vspace{-9mm}
\begin{center}
$^a$ Department of Physics, Urmia Branch, Islamic Azad University, Urmia, Iran\\
 \emph{*E-mail address}: alirezakhalili2002@yahoo.co.in (a.khalili@iaurmia.ac.ir)\\

$^b$ Department of Applied Mathematics and Theoretical Physics, University of Cambridge, Wilberforce Road, Cambridge, CB3 0WA, United Kingdom \\
 \emph{\dag E-mail address}: af454@cam.ac.uk\\
 $^c$ Department of Physics, University of Tabriz, Tabriz,  Iran\\
 \emph{\ddag E-mail address}: a.khalili@tabrizu.ac.ir\\
 
 $^{d}$  Department of Mathematics Cankaya University, 06530 Ankara, Turkey\\
$^{f}$  Institute of Space Sciences, P.O.BOX, MG-23, R 76900,~Magurele-Bucharest, Romania.\\ \emph{\S E-mail address}:
dumitru@cankaya.edu.tr

 \end{center}

\begin{abstract}

In this paper, we study $C^{\zeta}$-calculus on generalized Cantor sets, which have self-similar properties and fractional dimensions that exceed their topological dimensions. Functions with fractal support are not differentiable or integrable in terms of standard calculus, so we must involve local fractional derivatives. We have generalized the $C^{\zeta}$-calculus on the generalized Cantor sets known as middle-$\xi$ Cantor sets. We have suggested  a calculus on the middle-$\xi$ Cantor sets for different values of $\xi$ with $0<\xi<1$. Differential equations on the middle-$\xi$ Cantor sets have been solved, and we have presented the results using illustrative examples. The conditions for super-, normal, and sub-diffusion on fractal sets are given.

\end{abstract}




\section{Introduction}

It is well known that many phenomena in nature can be modeled by fractals; these shapes can be observed almost anywhere in the natural world \cite{Barnsley}. Fractal antennas have maximal  length, area, and volume to accommodate a multi-band or wide-band design, which is useful in cellular telephone and microwave communications \cite{Cohen,Cohen-2,Cohen-3}. Fractals also play important roles in biology. For example, in the neural and vascular networks of the human body, viruses and certain tumors grow and ramify in a fractal shape \cite{Gazit,Gazit-1,Gazit-2}. In these studies, researchers tried to predict and recognize osteoporosis from test results and from the fractal structure of bone texture \cite{Gazit-3}. Fractals have also been hypothesized to be important for human perception of beauty in artworks \cite{Taylor1999,Bountis2017}. Non-Markovian random walks and fractal dimensions which are connected to physical properties of fractal sets were studied in \cite{Balankin-1,Balankin-2,Alireza-cub}. The polynomial asymptotic behavior of the Wiener index on infinite lattices including fractals has been given \cite{Ori}.

Anomalous diffusion on fractals has received attention in recent years from various researchers \cite{Poirier,Ben}. Non-constant diffusion coefficients have been studied \cite{Petersen}, and the diffusion coefficient is proportional to a power of the noise intensity \cite{Robertwas}. The models for the diffusion coefficient characterize random motions in various regions of parameterized space and relate to the fractal structure as a function of the slope of the map \cite{Mark,Klages}. Fractal structural parameters have been applied to obtain porosity and tortuosity for micro-porous solids \cite{Schieferstein}. The scale-dependent fractal dimension for a random walk trajectory was used to derive the diffusion coefficient \cite{Gmachowski,Bujan}. The quenched-trap model on a fractal lattice does not lead to continuous-time random walks if the spectral dimension is less than 2 \cite{Miyaguchi}.

Fractional calculus has been applied to define derivatives on fractal curves \cite{book-100,jo-5,Kolwankar,Nigmatullin,Asad}. Fractional derivatives have non-local properties, so that they are used to model processes with memory effects \cite{Herrmann-2,Hilfer-2}.

Local fractional derivatives were suggested and applied from a physics perspective \cite{Gangal-1,Gangal-2} in a formalism called $C^{\zeta}$-calculus (or $C^{\zeta}$-C), which has also been generalized for unbounded and singular functions \cite{GolmankhanehGauge}. Schr\"{o}dinger equations on fractal curves were derived  using $C^{\zeta}$-C and Feynman path methods \cite{GolmankhanehSchr}. A mathematical model of diffraction was given for fractal sets \cite{Golmankhanehgrat}. Non-local derivatives were defined for fractal sets and applied in fractal mediums \cite{Golmankhaneh-107-k,Golmankhaneh-1087-k,Golmankhanehoo-3}. The Fokker-Planck equation for thick fractal absorbers was derived in view of $C^{\zeta}$-C \cite{Golmankhaneh-rr87-k}. Recently, as an application of the mathematical model, experimental and simulation results were utilized to model sub-diffusion and super-diffusion in physical processes \cite{Alexander255963}.

 By conducting research along these lines, we have generalized  $C^{\zeta}$-C to middle-$\xi$ Cantor sets.

The outline of the paper is as follows. In Section 2 we review $C^{\zeta}$-C and the basic tools required. In Section 3 we apply the $C^{\zeta}$-C on the middle-$\xi$ Cantor sets. In Section 4 we consider and solve some differential equations on middle-$\xi$ Cantor sets, and in Section 5 we consider diffusion processes on such sets. Section 4 is devoted to the conclusion.

\section{Basic tools in the fractal calculus  \label{int2}}
In this section, we review the Cantor-like sets and their properties \cite{Robert}, then summarize $C^{\zeta}$-C \cite{Gangal-1,Gangal-2} and non-local $C^{\zeta}$-C  \cite{Golmankhaneh-107-k,Golmankhaneh-1087-k}.

\subsection{Middle-$\xi$ Cantor sets}
Let us consider a unit interval $J=[0,1]$, and construct the middle-$\xi$ Cantor fractal set $C^{\xi}$ from it as follows.

In the first step, we remove an open interval of length $\xi$ from the exact middle of the interval $I$, to obtain:
\begin{equation}\label{mnu}
  C_{1}^{\xi}=\left[0,\frac{1-\xi}{2}\right]
  \cup\left[\frac{1+\xi}{2},1\right].
\end{equation}

In the second step, we pick up two open disjoint intervals with length $\xi^2$ from the middle of each of the remaining intervals which comprise the set $C_{1}^{\xi}$, to obtain
\begin{multline}\label{mnubv}
C_{2}^{\xi}=\left[0,\frac{1-\xi-2\xi^2}
{4}\right]\cup\left[\frac{1-\xi+2\xi^2}{4},
\frac{1-\xi}{2}\right] \\ \cup\left[\frac{1+\xi}{2},\frac{3+
\xi-2\xi^2}{4}\right]\cup\left[\frac{3+
\xi+2\xi^2}{4},1\right].
\end{multline}

After iterating this process infinitely many times, with the set constructed at stage $k$ being denoted by $C_{k}^{\xi}$, we obtain the definition of the middle-$\xi$ Cantor set as follows:
\begin{equation}\label{11mnu}
  C^{\xi}=\bigcap_{k=1}^{\infty}C_{k}^{\xi}.
\end{equation}
It is clear that the set $C^{\xi}$ has a self-similarity property, which makes it easy for us to find its fractional dimension. Namely, for every middle-$\xi$ Cantor set the  Hausdorff dimension is given by
\begin{equation}\label{xzaw}
  \dim_{H}(C^{\xi})=\frac{\log 2}{\log 2 -\log(1-\xi)},
\end{equation}
where $H(C^{\xi})$ is the  Hausdorff measure which was used to derive Hausdorff dimension \cite{Robert}. Furthermore, the middle-$\xi$ Cantor set has zero Lebesgue measure, because \cite{Robert}:
 \begin{equation}\label{xza}
  \textmd{L}_{m}(C^{\xi})=\lim_{k \rightarrow \infty}\textmd{L}_{m}(C_{k}^{\xi})=\lim_{k \rightarrow \infty}(1-\xi)^{k}=0.
\end{equation}

\begin{remark}
If we choose $\xi=1/3$, $\xi=1/4$, $\xi=1/5$, then we obtain the Cantor triadic set,  4-adic-type Cantor-like set, and 5-adic-type Cantor-like set, respectively. See Section 3 below for more details on these sets.
\end{remark}

\subsection{Local fractal calculus}
If $C^{\xi}$ is a middle-$\xi$ Cantor set contained in an interval $J=[v,w]\subset\Re$, then the flag function for $C^{\xi}$ is indicated by
$\varphi(C^{\xi},J)$ and  defined by \cite{Gangal-1,Gangal-2}
\begin{equation}
\varphi(C^{\xi},J)=\begin{cases}
    1 ~~~\textmd{if}~~~ C^{\xi}\cap J \neq \emptyset\\
    0 ~~~~\textmd{otherwise}.
\end{cases}
\end{equation}

For a set $C^{\xi}$ and a subdivision $Q_{[v,w]}=\{v=y_0,y_1,y_2,\dots,y_n=w\}$ of the interval $J=[v,w]$, we define
\begin{equation}
\rho^{\zeta}[C^{\xi},J]=\sum_{i=1}^{n}
\frac{(y_{i}-y_{i-1})^{\zeta}}{\Gamma
(\zeta+1)}~\varphi(C^{\xi},[y_{i-1},y_{i}])
\end{equation}
for any $\zeta$ with $0<\zeta\leq 1$. Given $\delta>0$, the associated coarse-grained mass function $\gamma_{\delta}^{\zeta}(C^{\xi},v,w)$ of the intersection $C^{\xi}\cap [v,w]$ is given by
\begin{equation}
\gamma_{\delta}^{\zeta}(C^{\xi},v,w)=\inf_{Q_{[v,w]}:
|Q|\leq\delta}\sigma^{\zeta}[C^{\xi},J],
\end{equation}
where the infimum is taken over all subdivisions $Q$ of $[v,w]$ satisfying $|Q|:=\max_{1\leq i\leq n}(y_{i}-y_{i-1})\leq\delta$. Then the \textbf{mass function} $\gamma^{\zeta}(C^{\xi},v,w)$ is given by \cite{Gangal-1,Gangal-2}
\begin{equation}\label{s}
    \gamma^{\zeta}(C^{\xi},v,w)=\lim_{\delta\rightarrow0} \gamma^{\zeta}_{\delta}(C^{\xi},v,w).
\end{equation}

The \textbf{integral staircase function} $S_{C^{\xi}}^{\zeta}(y)$ of order
$\zeta$ for a fractal set $C^{\xi}$ is defined in \cite{Gangal-1,Gangal-2} by
\begin{equation}\label{t}
    S_{C^{\xi}}^{\zeta}(y)=\begin{cases}
    \gamma^{\zeta}(C^{\xi},v_{0},y) ~~~~\text{if} ~~~~~y\geq a_{0}\\
    -\gamma^{\zeta}(C^{\xi},v_{0},y) ~~~~\text{otherwise},
\end{cases}
\end{equation}
where $v_{0}$ is an arbitrary real number. A point $y$ is a point of change of a function  $u(y)$ which is not constant over any open interval $(v,w)$ involving $y$. All points of change of $y$ is named the set of change
of $u(y)$ and is indicated by $\textbf{S}\textmd{ch}(S_{C^{\xi}}^{\zeta})$ \cite{Gangal-1,Gangal-2}. If $\textbf{S}\textmd{ch}(S_{C^{\xi}}^{\zeta})$ is a closed set and every point in it is a limit point, then $\textbf{S}\textmd{ch}(S_{C^{\xi}}^{\zeta})$ is called $\zeta$-perfect.

The $\varsigma$-dimension of $C^{\xi}\cap[v,w]$ is
\begin{align}\label{sa}
    \dim_{\varsigma}(C^{\xi}\cap[v,w])&=\inf\{\zeta:
    \gamma^{\zeta}(C^{\xi},v,w)=0\}\nonumber\\
    &=\sup\{\zeta:\gamma^{\zeta}(C^{\xi},v,w)=\infty\}.
\end{align}
For the definitions of $C^{\xi}$-limits and $C^{\xi}$-continuity, we refer the reader to \cite{Gangal-1,Gangal-2}.

\subsection{ $C^{\zeta}$-Differentiation}
If $C^{\xi}$ is an $\zeta$-perfect set, then the $C^{\zeta}$-derivative
of a function $u$ defined on $C^{\xi}$ at a point $y$ is defined to be the following, assuming the limit exists \cite{Gangal-1,Gangal-2}:
\begin{equation}
    D_{C^{\xi}}^{\zeta}u(y)=\begin{cases}
    C^{\xi}\text{-}\lim_{z\rightarrow y} \frac{f(z)-f(y)}{S_{C^{\xi}}^{\zeta}(z)-S_{C^{\xi}}^
    {\zeta}(y)}, ~~~\text{if }z\in C^{\xi},\\
    0, ~~~~~~~~~~~~~~~~~~~~~~~~~~~~~~~\textmd{otherwise}.
    \end{cases}
 \end{equation}

Let $u$ be a bounded function on $C^{\xi}$ and $J$ be a closed interval as above \cite{Gangal-1,Gangal-2}. Then we define
\begin{eqnarray}
  \mathfrak{M}[u,C^{\xi},J] &=& \sup_{y\in C^{\xi}\cap J} u(y)~~~~\text{if} ~~C^{\xi}\cap J\neq 0\\
  &=& 0 ~~~~~ \text{otherwise,}
\end{eqnarray}
and similarly
\begin{eqnarray}
  \mathfrak{m}[u,C^{\xi},J] &=& \inf_{y\in C^{\xi}\cap J} u(y)~~~~\text{if} ~~C^{\xi}\cap J\neq 0\\
  &=& 0 ~~~~~ \text{otherwise}
\end{eqnarray}

If $S_{C^{\xi}}^{\zeta}(y)$ is finite for $y\in [v,w]$, and $Q=\{v=y_0,y_1,\dots,y_n=w\}$ is a subdivision of $[v,w]$, then the upper $C^{\zeta}$-sum and lower $C^{\zeta}$-sum for a function $u$ over the subdivision $Q$ are given respectively by \cite{Gangal-1,Gangal-2}

\begin{equation}\label{y}
    \mathfrak{U}^{\zeta}[u,C^{\xi},Q]=\sum_{j=1}^{m}\mathfrak{M}[u,C^{\xi},
    [y_{j},y_{j-1}]](S_{C^{\xi}}^{\zeta}(y_{j})-
    S_{C^{\xi}}^{\zeta}(y_{j-1}))
\end{equation}
and
\begin{equation}\label{yo}
    \mathfrak{L}^{\zeta}[u,C^{\xi},Q]=\sum_{j=1}^{m}\mathfrak{m}[u,C^{\xi},
    [y_{j},y_{j-1}]]
    (S_{C^{\xi}}^{\zeta}(y_{j})-S_{C^{\xi}}^{\zeta}(y_{j-1}))
\end{equation}

If $u$ be a bounded function on $C^{\xi}$. we say that $u$ is
$C^{\zeta} $-integrable
 on $[v,w]$ if \cite{Gangal-1,Gangal-2} the two quantities
\begin{align}
\label{cft1}    \underline{\int_{v}^{w}}u(y) d_{C^{\xi}}^{\zeta}y&=\sup_{Q_{[v,w]}} \mathfrak{L}^{\zeta}[u,C^{\xi},Q], \\
\label{cft2}    \overline{\int_{v}^{w}}u(y) d_{C^{\xi}}^{\zeta}y&=\inf_{Q_{[v,w]}} \mathfrak{L}^{\zeta}[u,C^{\xi},Q],
\end{align}
are equal. In that case the $C^{\zeta}$-integral of $u$ on $[v,w]$ is  denoted
by $\int_{v}^{w}u(y) d_{C^{\xi}}^{\zeta}y$ and is given by the common value of \eqref{cft1},\eqref{cft2}.

\textit{Fundamental Theorem of $C^{\zeta}$-Calculus.} Suppose that $u(y):C^{\xi}\rightarrow \Re$   is $C^{\zeta}$-continuous  and  bounded  on $C^{\xi}$. If we define $\mathrm{g}(z)$ by
\begin{equation}\label{pokl}
  \mathrm{g}(z)=\int_{a}^{z}\mathrm{f}(y)d_{C^{\xi}}^{\zeta}y,
\end{equation}
for all $y\in [v,w]$, then:
\begin{equation}\label{desw}
  D_{C^{\xi}}^{\zeta}\mathrm{g}(y)=\mathrm{f}(y)\chi_{C^{\xi}}(y).
\end{equation}
where $\chi_{C^{\xi}}(y)$ is the characteristic function of the middle-$\xi$ Cantor set.

Conversely, if $\mathrm{f}(y)$ is an $C^{\zeta}$-differentiable function, then we have \cite{Gangal-1,Gangal-2}
\begin{equation}\label{desw3}
  D_{C^{\xi}}^{\zeta}\mathrm{f}(y)=\mathrm{h}(y)\chi_{C^{\xi}}(y).
\end{equation}
for some function $\mathrm{h}$, and consequently it follows that
\begin{equation}\label{bvc}
  \int_{v}^{w}\mathrm{h}(y)d_{C^{\xi}}^{\zeta}y=\mathrm{f}(b)-\mathrm{f}(a).
\end{equation}

\section{Staircase functions on middle-$\xi$ Cantor sets}
In this section, we plot middle-$\xi$ Cantor sets and their  staircase functions $S_{C^{\xi}}^{\zeta}(y)$ for special cases, in order to present details of the paper.

\subsection{The Cantor triadic set}

The Cantor triadic set is generated by iteration as follows:
\begin{itemize}
\item \textit{Step 1.} Remove an open interval of length $1/3$ from the middle of the interval $J=[0,1]$.
\item \textit{Step 2.} Remove an open interval of length $(1/3)^2$ from the middle of each one of the closed intervals with length $1/3$ remaining from step $1$.
\item \dots
\item \textit{Step k.} Remove an open interval of length $(1/3)^k$ from the middle of each one of the closed intervals with length $(1/3)^{k-1}$ remaining from step $k-1$.
\end{itemize}

In the case of the Cantor triadic sets, utilizing Eqs.(\ref{xzaw}) and (\ref{sa}), we get $\varsigma$-dimension as follows:
\begin{equation}\label{vfx6}
  \dim_{\varsigma}(C^{1/3}\cap[v,w])=\dim_{H}(C^{1/3})=0.63.
\end{equation}

\begin{figure}	
	\centering
	\begin{subfigure}[t]{0.4\textwidth}
		\centering
		\includegraphics[width=\textwidth]{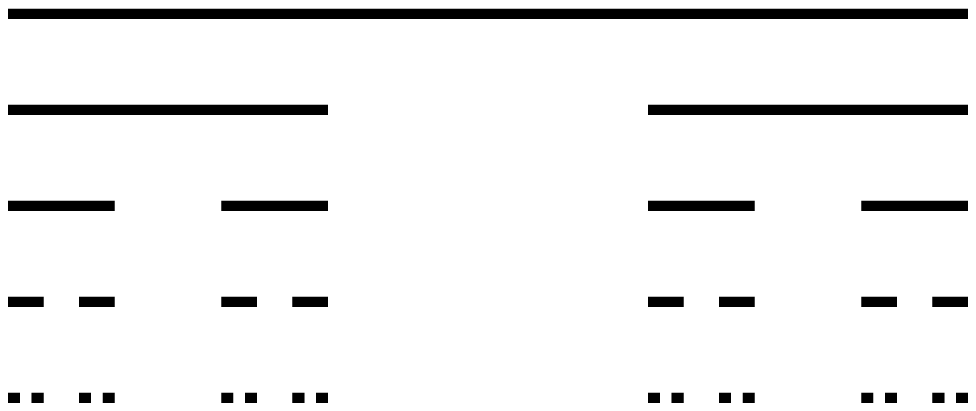}
		\caption{Steps 0--4 of the generating process for the Cantor triadic set}\label{rf1:a}
	\end{subfigure}
	\quad
	\begin{subfigure}[t]{0.4\textwidth}
		\centering
		\includegraphics[width=\textwidth]{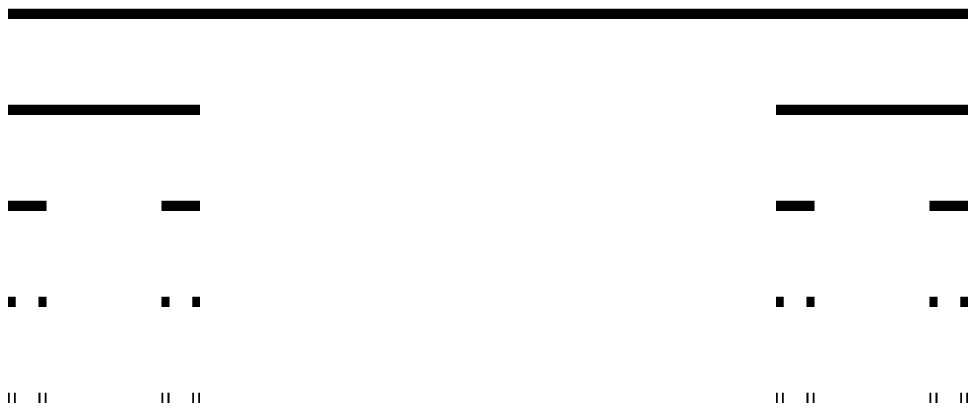}
		\caption{Steps 0--4 of the generating process for the 5-adic-type Cantor-like set}\label{rf1:b}
	\end{subfigure}
	\\
	\begin{subfigure}[b]{0.4\textwidth}
		\centering
		\includegraphics[width=\textwidth]{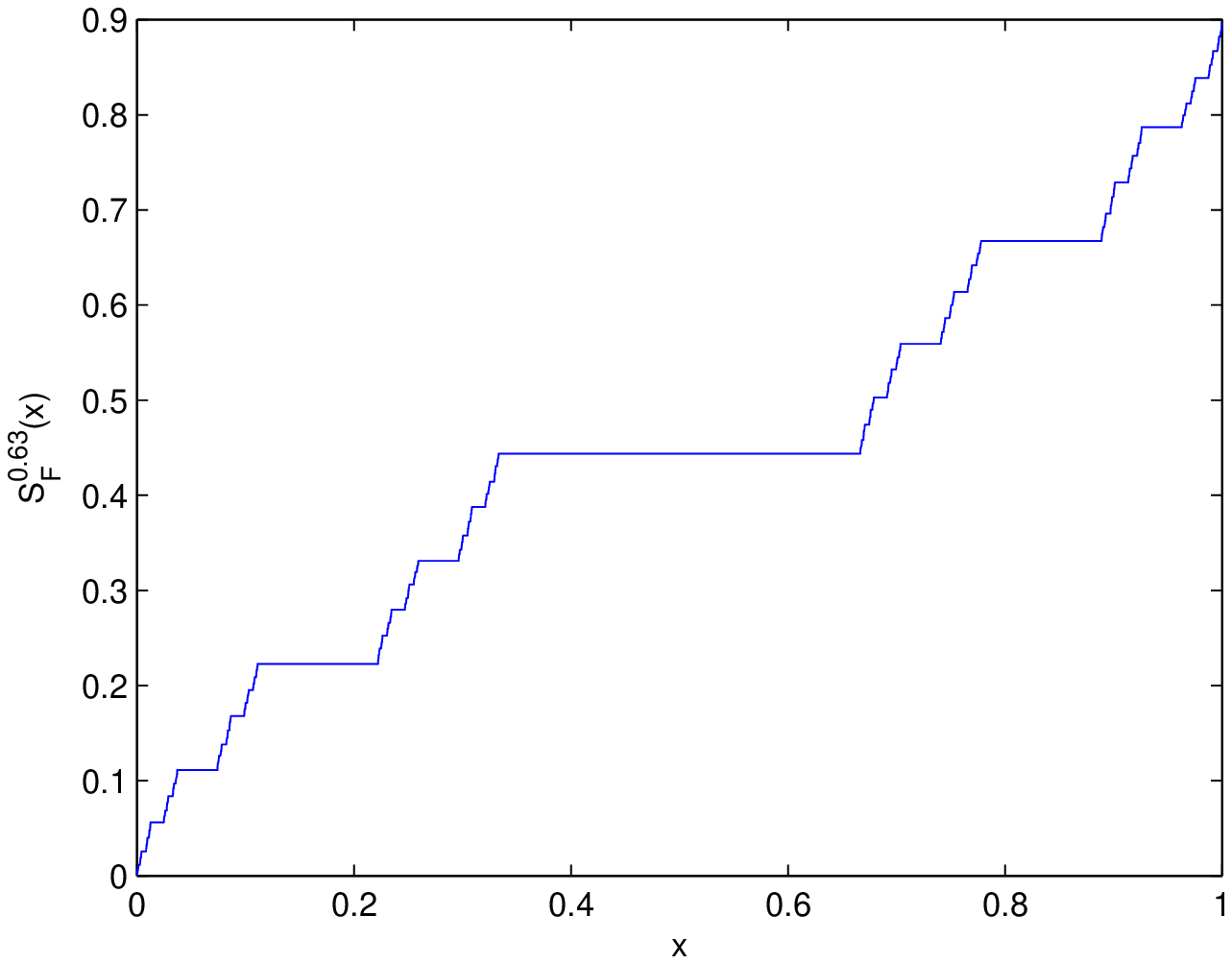}
		\caption{The staircase function corresponding to the Cantor triadic set}\label{rf1:c}
	\end{subfigure}
	\quad
	\begin{subfigure}[b]{0.4\textwidth}
		\centering
		\includegraphics[width=\textwidth]{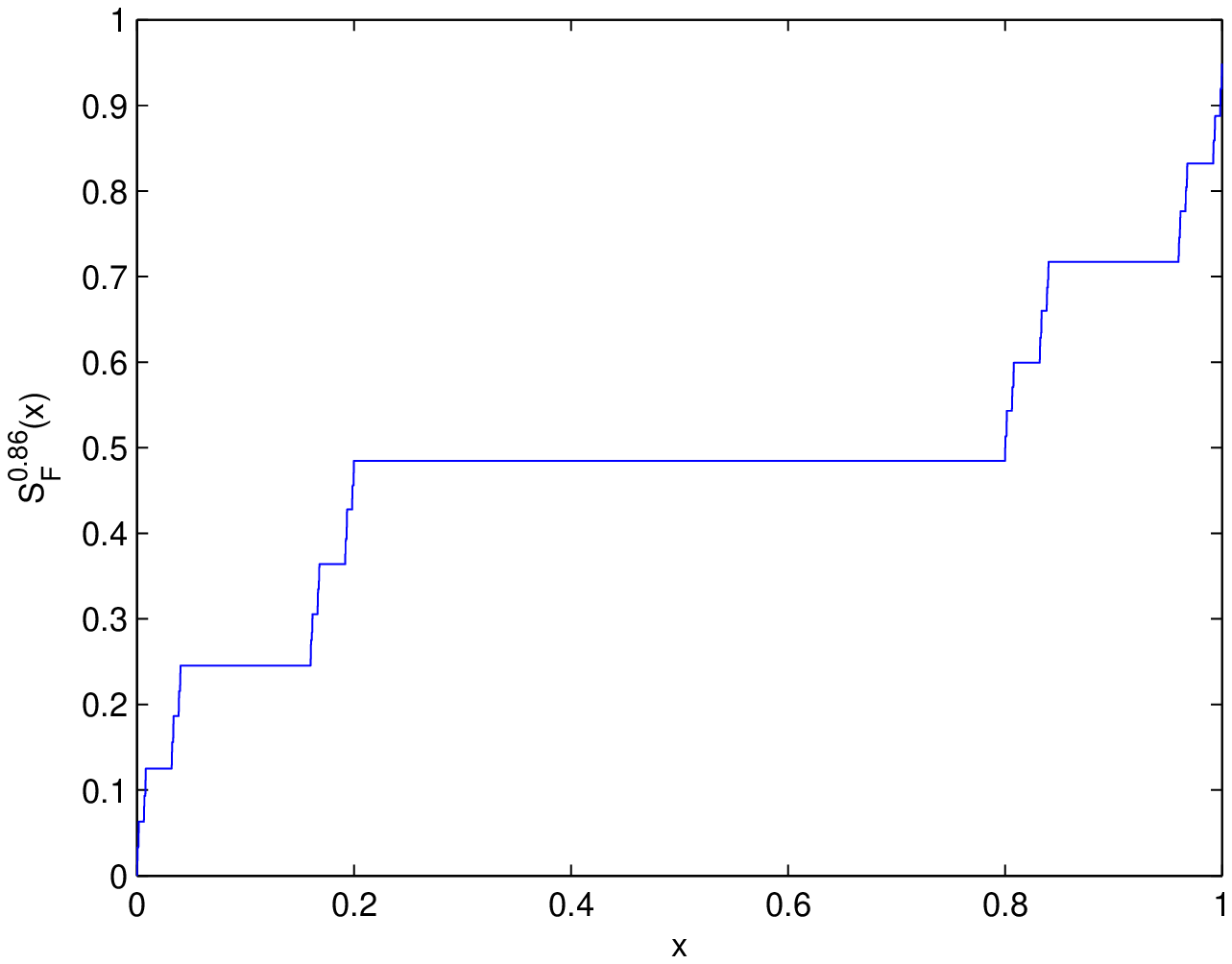}
		\caption{The staircase function corresponding to the 5-adic-type Cantor-like set}\label{rf1:d}
	\end{subfigure}
	\caption{Basic properties of some example Cantor sets}\label{rf1}
\end{figure}

In Figure \ref{rf1:a} we draw the process mentioned above which established the Cantor triadic set.

Using  Eq. (\ref{t}) the staircase function of  the Cantor triadic set $(S_{F^{1/3}}^{0.63}(x))$ is sketched in Figure \ref{rf1:c}.

The staircase functions have important roles in $C^{\zeta}$-C, being used in integration and differentiation of functions with fractal support.

\subsection{The 5-adic-type Cantor-like set}
The procedure to achieve the 5-adic-type Cantor-like set is similar to that for the Cantor triadic set, only differed to remove the open interval of length $1/5$ in every stage. We exhibit these steps in Figure 1b.

The $\varsigma$-dimension of the  5-adic-type Cantor-like set, considering Eqs. (\ref{xzaw}) and (\ref{sa}), is
\begin{equation}\label{vfx8}
  \dim_{\varsigma}(C^{1/5}\cap[v,w])=\dim_{H}(C^{1/5})=0.86.
\end{equation}
From  Eq. (\ref{t}), the staircase function of the 5-adic-type Cantor-like $(S_{C^{1/5}}^{0.86}(y))$ is plotted in Figure 1d.

\section{Differential equations on middle-$\xi$ Cantor sets}
In this section, first, we study the integration and differentiation of functions whose support is a middle-$\xi$ Cantor sets. Secondly, differential equations formulated on middle-$\xi$ Cantor sets are suggested and solved using illustrative examples.

\begin{example} \label{S3:ex1}
Consider a function with the fractal Cantor triadic set support as follows:
\begin{equation}\label{x23}
  f(x)=\sin\left((2\pi x\chi_{C^{1/3}}^{0.63}(x)\right)
\end{equation}
where $\chi_{C^{1/3}}^{0.63}$ is the characteristic function of the fractal Cantor triadic set. We plot the function $f(x)$ in Figure 2a.

\begin{figure}	
	\centering
	\begin{subfigure}[t]{0.4\textwidth}
		\centering
		\includegraphics[width=\textwidth]{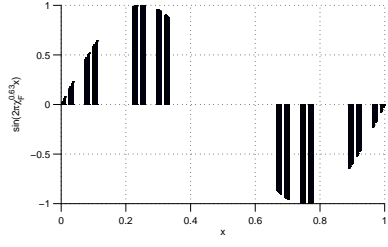}
		\caption{Graph of the function $f(x)$}\label{frtswe1:a}
	\end{subfigure}
	\\
	\begin{subfigure}[b]{0.4\textwidth}
		\centering
		\includegraphics[width=\textwidth]{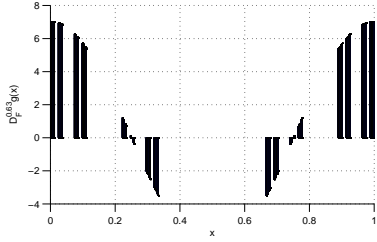}
		\caption{Graph of the fractal derivative $D_{F^{1/3}}^{0.63}f(x)$}\label{frtswe1:b}
	\end{subfigure}
	\quad
	\begin{subfigure}[b]{0.4\textwidth}
		\centering
		\includegraphics[width=\textwidth]{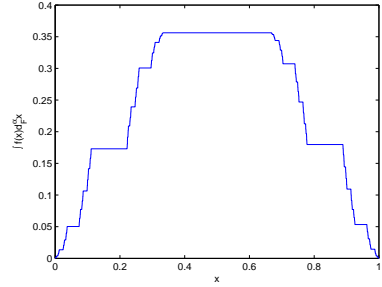}
		\caption{Graph of the fractal integral $\int_{0}^{x} f(x')d_{F^{1/3}}^{0.63}x'$}\label{frtswe1:c}
	\end{subfigure}
	\caption{Graphs relevant to Example \ref{S3:ex1}}\label{frtswe1}
\end{figure}


The $C^{\zeta}$-derivative of $f(x)$, using conjugacy of $C^{\zeta}$-C and ordinary calculus \cite{Gangal-1,Gangal-2}, is as follows:
\begin{equation}\label{x23}
  D_{C^{1/3}}^{0.63}f(x)=2\Gamma(1+0.63)\pi \cos\left(2\pi x \chi_{C^{1/3}}^{0.63}(x)\right),
\end{equation}
where $\Gamma(\cdot)$ denotes the gamma function. A plot of the function $D_{C^{1/3}}^{0.63}f(x)$ is shown  in Figure 2b.

The $C^{\zeta}$-integration of $f(x)$, considering conjugacy of $C^{\zeta}$-C between the ordinary calculus \cite{Gangal-1,Gangal-2}, will be as follows:
\begin{align*}
\int_{0}^{1} \sin\left(2\pi x \chi_{C^{1/3}}^{0.63}(x)\right)d_{C^{1/3}}^{0.63}x&=\frac{-1}{2\pi \Gamma(1+0.63)} \left[\cos(2\pi\Gamma(1+0.63) S_{C^{1/3}}^{0.63}(x))\right]_{0}^{1} \\ &=\frac{1}{2\pi \Gamma(1+0.63)}\left[1 - \cos(2\pi\Gamma(1+0.63) S_{C^{1/3}}^{0.63}(1))\right] \\
&=0,
\end{align*}
where we use $\Gamma(1+0.63) S_{C^{1/3}}^{0.63}(1)=1$. In Figure 2c we plot the integral function of $f(x)$ over $[0,1]$.

\end{example}

\begin{example} \label{S3:ex2}
Suppose we have a function on the fractal 5-adic-type Cantor-like set as follows:
\begin{equation}\label{x225873}
  g(x)= x^2\chi_{C^{1/5}}^{0.86}(x).
\end{equation}
This function  $g(x)$ is sketched  in Figure 3a.

The $C^{\xi}$-derivative of $g(x)$ is derived by a similar method as used in Example \ref{S3:ex1}, which yields the following result:
\begin{equation}\label{x23de14}
  D_{C^{1/5}}^{0.86}g(x)=2x\chi_{C^{1/5}}^{0.86}(x).
\end{equation}
We plot $D_{F^{1/5}}^{0.86}g(x)$ in Figure 3b.

\begin{figure}	
	\centering
	\begin{subfigure}[t]{0.4\textwidth}
		\centering
		\includegraphics[width=\textwidth]{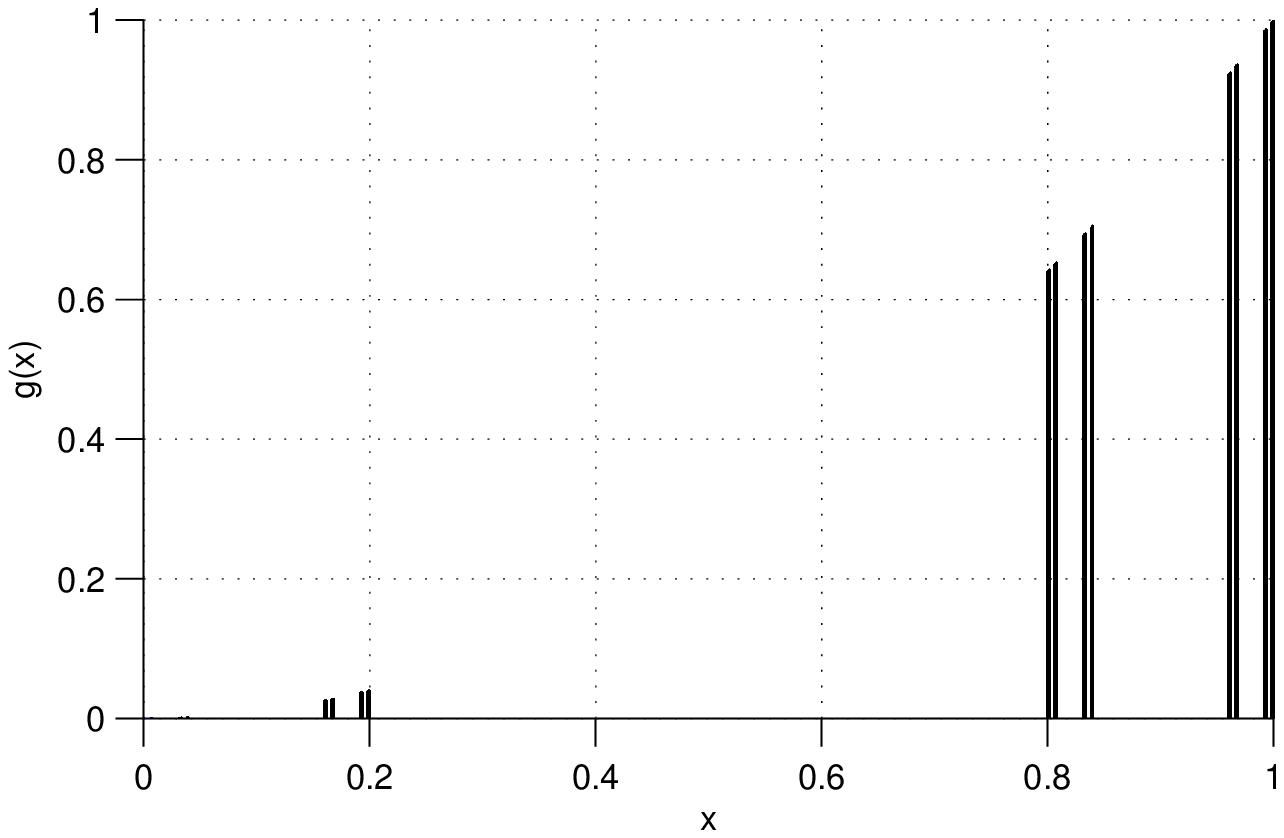}
		\caption{Graph of the function $g(x)$}\label{1gytderq:a}
	\end{subfigure}
	\\
	\begin{subfigure}[b]{0.4\textwidth}
		\centering
		\includegraphics[width=\textwidth]{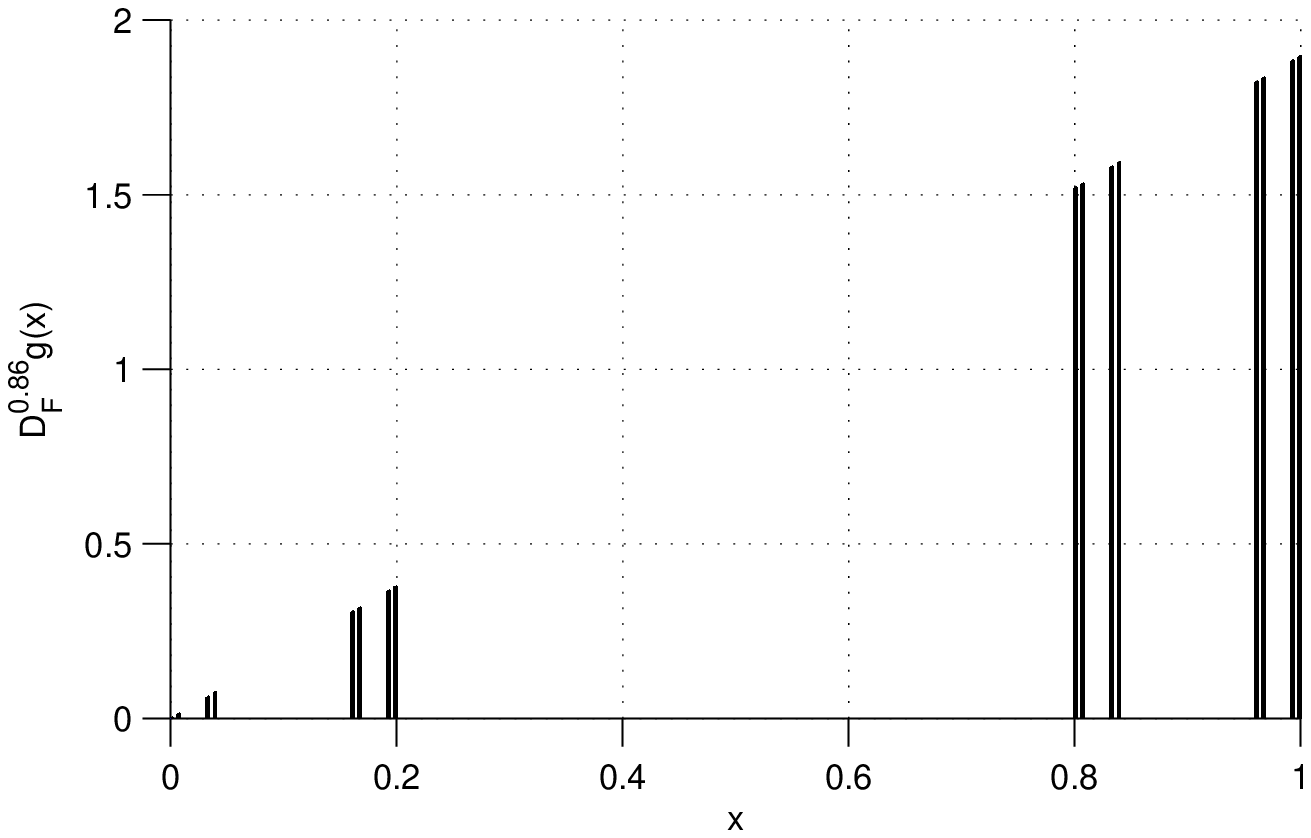}
		\caption{Graph of the fractal derivative $D_{F^{1/5}}^{0.86}g(x)$}\label{1gytderq:b}
	\end{subfigure}
	\quad
	\begin{subfigure}[b]{0.4\textwidth}
		\centering
		\includegraphics[width=\textwidth]{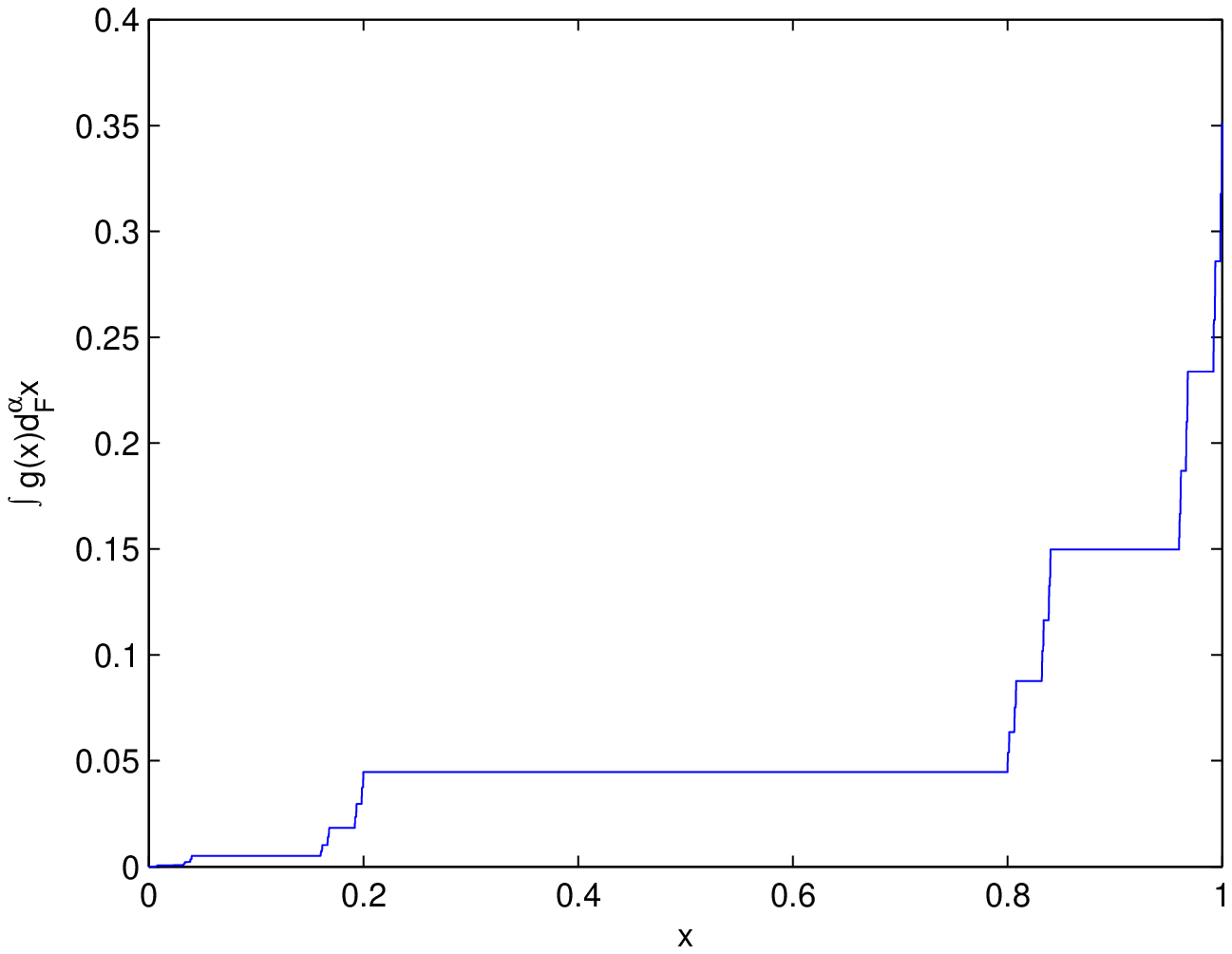}
		\caption{Graph of the fractal integral $\int_{0}^{x} g(x')d_{F^{1/5}}^{0.86}x'$}\label{1gytderq:c}
	\end{subfigure}
	\caption{Graphs relevant to Example \ref{S3:ex2}}\label{1gytderq}
\end{figure}


In the same manner, we obtain the $F^{\zeta}$-integral of $g(x)$ as follows:
\begin{eqnarray}\label{xgtfy23}
  \int_{0}^{1} (\chi_{C^{1/5}}^{0.86})x^2 d_{C^{1/5}}^{0.86}x&=& \frac{1}{3}\left(S_{C^{1/5}}^{0.86}(x)\right)^3\bigg|_{0}^{1}\nonumber\\
  &=&\frac{1}{3}\left(S_{C^{1/5}}^{0.86}(1)\right)^3=0.2846
\end{eqnarray}
We plot the $C^{\zeta}$-integral of $g(x)$ in Figure 3c.
\end{example}
We are going to use the above results and examples for solving differential equations on the middle-$\xi$ Cantor sets.

\section{Diffusion on middle-$\xi$ Cantor Sets }
In this section, we define and give conditions for super, normal and sub-diffusion on middle-$\zeta$ Cantor sets.

\subsection{Super-diffusion}
Let us consider time as continuous and space as a middle-$\xi$ Cantor set. We consider a probability function $W(x,t)$ which is an $C^{\zeta}$-differentiable function of the space coordinate $x$ and is a differentiable function of time $t$ in the sense of standard calculus. The fractal diffusion equation for a random walk is suggested with the conditional probability $W(x,t)$ as follows:
\begin{equation}\label{mex1}
\chi _{C^{\xi}}(x)~\frac{\partial W(x,t)}{\partial t}= K_{C^{\zeta}} (D_{C^{\xi},x}^{\zeta})^2 W(x,t),~~~t\in \Re,~x \in C^{\xi},
\end{equation}
with the initial condition
\begin{equation}\label{de}
  W(x,t=0)=\delta_{C^{\xi}}^{\zeta}(x),
\end{equation}
where $[K_{C^{\xi}}]=(Length^{2\zeta}/Time)$ is a generalized diffusion coefficient, $\delta_{C^{\xi}}^{\zeta}(x)$ is the Dirac delta function with fractal support. Using conjugacy of $C^{\zeta}$-C between standard calculus \cite{Gangal-1,Gangal-2}, we have the solution for Eq. (\ref{mex1}) as follows:
\begin{equation}\label{dsmex}
  W(x,t)=\frac{t^{-1/2}}{\sqrt{4 \pi K_{C^{\xi}}}}\exp\left[\frac{-S_{C^{\xi}}^{\zeta}(x)^2}{4 K_{C^{\xi}} t}\right].
\end{equation}
Since $S_{C^{\xi}}^{\zeta}(x)\leq x^{\zeta}$, then Eq. (\ref{dsmex}) can be written as:
\begin{equation}\label{zswe}
 W(x,t)\mapsto\frac{t^{-1/2}}{\sqrt{4 \pi K_{C^{\xi}}}}  \exp\left[\frac{-x^{2\zeta}}{4 K_{C^{\xi}}t} \right].
\end{equation}
The function $W(x,t)$ is indicated as the probability distribution of super-diffusion on Cantor sets.

Accordingly, the mean square random walk is
\begin{equation}\label{aws2uy}
  \langle S_{C^{\xi}}^{\zeta}(x)^2\rangle=4 K_{C^{\xi}} t.
\end{equation}
Using the upper bound $S_{C^{\xi}}^{\zeta}(x)^2\leq x^{2\zeta}$, we have
\begin{equation}\label{aws2uy}
 \langle x^2\rangle \mapsto 4 K_{C^{\xi}} t^{1/\zeta}.
\end{equation}

\subsection{Normal diffusion}
Let us consider space as a middle-$\xi$ Cantor set and fractal time associated with the middle-$\xi$ Cantor set, both sets having the same value of $\xi$ and the same dimension $\zeta$. The random walk conditional probability $W(x,t)$ is given by
\begin{equation}\label{bvgtf}
 D_{C^{\xi},t}^{\zeta} W(x,t)= G_{C^{\xi}} (D_{C^{\xi},x}^{\zeta})^2 W(x,t),
\end{equation}
where $[G_{C^{\xi}}]= (Length^{2\zeta}/Time^{\zeta})$ is a diffusion coefficient. The solution for Eq. (\ref{bvgtf}) with the initial condition Eq. (\ref{de}), utilizing conjugacy of $C^{\zeta}$-C between standard calculus, is:
\begin{equation}\label{78pki}
  W(x,t)=\frac{S_{C^{\xi}}^{\zeta}(t)^{-1/2}}{\sqrt{4 \pi G_{C^{\xi}}}}\exp\left[\frac{-S_{C^{\xi}}^{\zeta}(x)^2}{4 G_{C^{\xi}} S_{C^{\xi}}^{\zeta}(t)}\right]
\end{equation}
Considering the upper bound on the $S_{C^{\xi}}^{\zeta}(\cdot)$, we obtain:
\begin{equation}\label{x9}
   W(x,t)\mapsto \frac{t^{-\zeta/2}}{\sqrt{4 \pi G_{C^{\xi}}}}\exp\left[\frac{-x^{2 \zeta}}{4G_{C^{\xi}} t^{\zeta}}\right]
\end{equation}
The function $W(x,t)$ indicates the probability distribution for normal diffusion with a non-Gaussian propagator.

Applying conjugacy of $C^{\zeta}$-C between standard calculus, we arrive at the mean square of displacement
\begin{equation}\label{ew96}
  \langle S_{C^{\xi}}^{\zeta}(x)^2\rangle=4 G_{C^{\xi}} S_{C^{\xi}}^{\zeta}(t),
\end{equation}
and if we use the upper bound on $ S_{C^{\xi}}^{\zeta}(\cdot)$ we can write Eq. (\ref{ew96}) as follows:
\begin{equation}\label{ty}
\langle x^2\rangle\mapsto4 G_{C^{\xi}} t.
\end{equation}

\subsection{Sub-diffusion}
Let us consider time as a middle-$\xi$ Cantor set with dimension $\beta$ and space as a middle-$\xi$ Cantor set with dimension $\zeta$. A random walk on this fractal space-time has conditional probability that can be obtained by the following differential equation:
\begin{equation}\label{o}
   D_{C^{\xi'},t}^{\beta} W(x,t)= \chi_{C^{\xi'}} L_{C^{\xi}} (D_{C^{\xi},x}^{\zeta})^2 W(x,t)
\end{equation}
where $[L_{C^{\xi}}]=(Length^{2\zeta}/Time^{\beta})$ is a diffusion coefficient. Solving  Eq. (\ref{bvgtf}) with the initial condition Eq. (\ref{de}), using  conjugacy of $C^{\zeta}$-C between standard calculus, one can obtain
\begin{equation}\label{derf}
  W(x,t)=\frac{S_{C^{\xi}}^{\beta}(t)^{-1/2}}{\sqrt{4 \pi L_{C^{\xi}}}}\exp\left[\frac{-S_{C^{\xi}}^{\zeta}(x)^2}{4 L_{C^{\xi}} S_{C^{\xi}}^{\beta}(t)}\right].
\end{equation}
In view of the upper bounds on $ S_{C^{\xi}}^{\zeta}(\cdot)$, we get
\begin{equation}\label{ewsq}
  W(x,t)\mapsto \frac{t^{-\beta/2}}{\sqrt{4 \pi L_{C^{\xi}}}}\exp\left[\frac{-x^{2\zeta}}{4 L_{C^{\xi}} t^{\beta}}\right].
\end{equation}
The function $W(x,t)$ is named as the probability of sub-diffusion for a random walk as indicated above.

Similarly to the previous cases, the mean square of displacement in this case will be
\begin{equation}\label{eewwqw96}
  \langle S_{C^{\xi}}^{\zeta}(x)^2\rangle=4 L_{C^{\xi}} S_{C^{\xi}}^{\beta}(t),
\end{equation}
and in the same manner we use upper bounds on $ S_{C^{\xi}}^{\zeta}(\cdot)$ to get
\begin{equation}\label{tytgvfrt}
\langle x^2\rangle\mapsto4 L_{C^{\xi}} t^{\beta/\zeta}.
\end{equation}

\begin{example}
Consider a random walk model on the fractal 5-adic-type Cantor-like set. The corresponding mean square value displacement of the random walk is given by:
\begin{equation}\label{lopm}
\langle S_{C^{\xi}}^{0.63}(x)^2\rangle=4 L_{C^{\xi}} S_{C^{\xi}}^{\beta}(t),
\end{equation}
or
\begin{equation}\label{tlomnuytgvfrt}
\langle x^2\rangle\mapsto4 L_{C^{\xi}} t^{\beta/0.86}.
\end{equation}
where the respective cases $\beta>0.86$, $\beta<0.86$, and $\beta=0.86$ are called super-diffusion,  sub-diffusion and normal diffusion on the fractal 5-adic-type Cantor-like set, respectively.

\begin{figure}[H]
  \centering
  \includegraphics[scale=0.4]{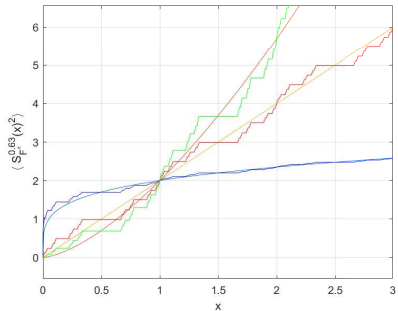}
 \caption{Graph of the mean value displacement  of the random walk model supper-, sub-diffusion and normal diffusion in the case $\alpha=0.86$  }\label{hbvcxs7458514}
\end{figure}

In Figure \ref{hbvcxs7458514} we draw the graphs of mean square value displacement random walk model for super-diffusion, sub-diffusion, and normal diffusion on fractal 5-adic-type Cantor-like sets in the case $\zeta=0.86$.
\end{example}

\begin{remark}
We conclude the results of Section 5  as follows.
\begin{enumerate}
\item The diffusion is super-diffusion on the middle-$\xi$ Cantor set if $\zeta<\beta$.
\item The diffusion is normal on the middle-$\xi$ Cantor set if $\zeta=\beta$.
\item The diffusion is sub-diffusion on the middle-$\xi$ Cantor set if $\zeta>\beta$.
\end{enumerate}
\end{remark}
\begin{remark}
In some figures we have plotted bars instead of points for the graphs of functions with fractal support, in order to make the results more clear.
\end{remark}

\section{Conclusion}
The $C^{\zeta}$-calculus is a generalization of ordinary calculus that can be applied on middle-$\xi$ Cantor sets for different values of $\xi$. Functions with middle-$\xi$ Cantor set support were considered, and their derivatives and integrals were derived  using $C^{\zeta}$-calculus, which shows the advantage of using $C^{\zeta}$-calculus over standard calculus. $C^{\zeta}$-derivatives on new fractal sets were discussed and compared for functions with different fractal supports. New differential equations involving  $C^{\zeta}$-derivatives on middle-$\xi$ Cantor sets have been suggested, that can be used as mathematical models for many physical problems. For example, we suggest conditions for super-, normal, and sub-diffusion on fractal sets.


\vspace{6pt}

\end{document}